\title{Hypermap operations of finite order}

\author{
Gareth A. Jones\\
School of Mathematics\\
University of Southampton\\
Southampton SO17  1BJ\\
U.K.\\
{\tt G.A.Jones@maths.soton.ac.uk}
\and
Daniel Pinto\\
Departamento de Matem\'atica\\
Universidade de Coimbra\\
Apartado 3008\\
3001--454 Coimbra\\
Portugal\\
{\tt dpinto@mat.uc.pt}
}

\documentclass[12pt]{article}
\usepackage[latin1]{inputenc}
\usepackage{a4wide}
\usepackage{latexsym}
\usepackage{amsmath}
\usepackage{amssymb}
\usepackage{amsfonts}
\newtheorem{thm}{Theorem}[section]

\newtheorem{cor}[thm]{Corollary}
\newtheorem{prop}[thm]{Proposition}
\date{}
\begin{document} 

\maketitle

\begin{abstract}
\noindent Duality and chirality are examples of operations of order $2$ on hypermaps. James showed that the groups of all operations on hypermaps and on oriented hypermaps can be identified with the outer automorphism groups ${\rm Out}\,\Delta\cong PGL_2({\bf Z})$ and ${\rm Out}\,\Delta^+\cong GL_2({\bf Z})$ of the groups $\Delta = C_2*C_2*C_2$ and $\Delta^+ = F_2$. We will consider the elements of finite order in these two groups, and the operations they induce.
\end{abstract}



\noindent{\bf Running head:} Hypermap operations\\

\section{Introduction}

A hypermap can be regarded as a transitive permutation representation $\Delta\to{\rm Sym}\,\Phi$ of the group
\[\Delta=\langle r_0, r_1, r_2\mid r_0^2=r_1^2=r_2^2=1\rangle\cong C_2*C_2*C_2,\]
the free product of three cyclic groups of order $2$, on a set $\Phi$ representing its flags; similarly an oriented hypermap (without boundary) can be regarded as a transitive permutation representation of the subgroup
\[\Delta^+=\langle \rho_0, \rho_1, \rho_2\mid \rho_0\rho_1\rho_2
=1\rangle=\langle \rho_0,\rho_2\mid -\rangle\cong F_2\]
of index $2$ in $\Delta$ (a free group of rank $2$) consisting of the elements of even word-length in the generators $r_i$, where $\rho_0=r_1r_2$, $\rho_1=r_2r_0$ and $\rho_2=r_0r_1$. In the first case the hypervertices, hyperedges and hyperfaces ($i$-dimensional constituents for $i=0, 1, 2$) are the orbits of the dihedral subgroups $\langle r_1, r_2\rangle$, $\langle r_2, r_0\rangle$ and $\langle r_0, r_1\rangle$, and in the second case they are the orbits of the cyclic subgroups $\langle\rho_0\rangle$, $\langle\rho_1\rangle$ and $\langle\rho_2\rangle$, with incidence given by nonempty intersection in each case. The local orientation around each hypervertex, hyperedge or hyperface is given by the cyclic order within the corresponding cycle of $\rho_0, \rho_1$ or $\rho_2$.

If $\cal H$ is a hypermap corresponding to a permutation representation $\theta:\Delta\to{\rm Sym}\,\Phi$, and if $\alpha$ is an automorphism of $\Delta$, then $\alpha^{-1}\circ\theta:\Delta\to{\rm Sym}\,\Phi$ corresponds to a hypermap ${\cal H}^{\alpha}$. (We invert $\alpha$ so that the automorphism $\alpha\circ\beta$, composed from left to right, sends $\cal H$ to $({\cal H}^{\alpha})^{\beta}$ and not $({\cal H}^{\beta})^{\alpha}$.) The hypervertices of ${\cal H}^{\alpha}$ are therefore the orbits of $\langle r_1^{\alpha}, r_2^{\alpha}\rangle$ on $\Omega$, and similarly for the hyperedges and hyperfaces. If $\alpha$ is an inner automorphism then ${\cal H}^{\alpha}\cong{\cal H}$ for all $\cal H$, so we have an induced action of the outer automorphism group ${\rm Out}\,\Delta = {\rm Aut}\,\Delta/{\rm Inn}\,\Delta$ as a group $\Omega$ of operations on isomorphism classes of hypermaps. James~\cite{Jam} has shown that this action is faithful, even when restricted by finite hypermaps. The same remarks apply to oriented hypermaps, with ${\rm Out}\,\Delta^+$ acting as a group $\Omega^+$ of operations. In either case, we will let $\omega_{\alpha}$ denote the operation on hypermaps (or oriented hypermaps) induced by an automorphism $\alpha$, or more precisely by the outer automorphism class $[\alpha]$ containing $\alpha$.

\medskip

\noindent{\bf Example 1.} For each $2$-element subset $\{i, j\}\subset\{i, j, k\}=\{0, 1, 2\}$, let $\alpha$ be the automorphism $\alpha_{ij}=\alpha_{ji}$ of $\Delta$ which transposes $r_i$ and $r_j$ and fixes $r_k$. Then $\omega_{\alpha}$ is the duality operation $\omega_{ij}=\omega_{ji}$ on hypermaps which transposes their $i$- and $j$-dimensional constituents while preserving their $k$-dimensional constituents. More generally, the symmetric group $S_3$ acts as a group of automorphisms of $\Delta$ by permuting the generators $r_i$, giving a faithful action of $S_3$ on hypermaps; the images of a hypermap under this action are called its {\sl associates}. Note that because of the inversion mentioned in the preceding paragraph, if $\alpha$ corresponds to the  $3$-cycle $(0,1,2)$, acting on subscripts, then the hypervertices of ${\cal H}^{\alpha}$ are the hyperedges of ${\cal H}$, and so on. An action of $S_3$ on oriented hypermaps has been described by Mach\`\i{} in~\cite{Mac}; see \S 6 for the relationship between these two actions (which is not entirely straightforward). 

\medskip

\noindent{\bf Example 2.} Since $\Delta^+$ is a normal subgroup of $\Delta$, conjugation by $r_1$ induces an automorphism
\[\alpha_r:\rho_2\mapsto\rho_2^{-1},\;\rho_0\mapsto\rho_0^{-1}\]
of $\Delta^+$. The corresponding operation $\omega_r$ reflects each oriented hypermap, reversing its orientation. (It is, in fact, induced by conjugation by any element of $\Delta\setminus\Delta^+$, since any two of them differ by an element of $\Delta^+$.) An oriented hypermap is called {\sl reflexible\/} or {\sl chiral\/} as it is or is not invariant under $\omega_r$.
 
\medskip

\noindent{\bf Example 3.} It was shown by Hall~\cite{Hal} that the free group $F_2$ has $19$ normal subgroups with quotient group $A_5$, corresponding to the $19$ orbits of ${\rm Aut}\,A_5=S_5$ on generating pairs for $A_5$. It follows that there are $19$ orientably regular  hypermaps with orientation-preserving automorphism group $A_5$, and these are described in~\cite{BJ}. Neumann and Neumann~\cite{NN} have shown that ${\rm Aut}\,F_2$ acts as the group $S_9\times S_{10}$ on these $19$ normal subgroups, with two orbits of lengths $9$ and $10$, so $\Omega^+={\rm Out}\,\Delta^+$ acts in the same way as a group of operations, with orbits ${\cal O}_1$ and ${\cal O}_2$ of lengths $9$ and $10$ on the corresponding hypermaps. Direct calculation (see Example~4 in \S 4) shows that ${\cal O}_1$ consists of three sets of three associates; their types are the cyclic permutations of $(5,2,5)$ (corresponding to the great dodecahedron of genus $4$), $(5,3,3)$ of genus $5$ and $(5,5,3)$ of genus $9$. The orbit ${\cal O}_2$ consists of three sets of associates of size $6$, $3$ and $1$, and their types are the cyclic permutations of $(5,2,3)$ (the isosahedron of genus $0$), $(5,3,3)$ of genus $9$ and $(5,5,5)$ of genus $13$. The hypermaps of type $(5,5,3)$ in the orbits ${\cal O}_1$ and ${\cal O}_2$ represent $\rho_0$ and $\rho_1$ by $5$-cycles which are conjugate in $S_5$ by odd and even permutations respectively. See also~\cite{HM, Pro} for generalisations from $A_5$ to various other finite groups.

\section{Automorphism groups of $\Delta$ and $\Delta^+$}

It is well known (see~\cite[Ch.~I, Prop.~4.5]{LS} for instance) that  ${\rm Out}\,F_2\cong GL_2({\bf Z})$, this group being represented faithfully on the abelianised group $F_2^{\rm ab}\cong {\bf Z}^2$, so $\Omega^+={\rm Out}\,\Delta^+\cong GL_2({\bf Z})$. We will choose this isomorphism so that $GL_2({\bf Z})$ acts on row vectors with respect to the ordered basis of $(\Delta^+)^{\rm ab}\cong {\bf Z}^2$ induced by the generators $\rho_2$ and $\rho_0$ of $\Delta^+$. Similarly, James has shown in~\cite{Jam} that  $\Omega={\rm Out}\,\Delta\cong PGL_2({\bf Z})=GL_2({\bf Z})/\{\pm I\}$. In order to understand the relationship between these two outer automorphism groups, and to find their periodic elements, we will use the facts that $GL_2({\bf Z})$ is generated by the matrices
\[X=\Bigl(\,\begin{matrix}0&-1\\ 1&0\end{matrix}\,\Bigr),
\quad Y=\Bigl(\,\begin{matrix}0&-1\\ 1&1\end{matrix}\,\Bigr)
\quad{\rm and}\quad
T=\Bigr(\,\begin{matrix}0&1\\ 1&0\end{matrix}\,\Bigr)\]
with defining relations
\[X^2=Y^3,\quad X^4=T^2=(XT)^2=(YT)^2=I\eqno(2.1)\]
(see~\cite[\S 7.2]{CM}), and that $PGL_2({\bf Z})$ is generated by the images $x=\pm X,\, y=\pm Y$ and $t=\pm T$ of $X, Y$ and $T$ with defining relations
\[x^2=y^3=t^2=(xt)^2=(yt)^2=1.\eqno(2.2)\]

Since $\Delta^+$ is the unique torsion-free subgroup of index $2$ in $\Delta$, it is a characteristic subgroup, and hence each automorphism $\alpha$ of $\Delta$ restricts to an automorphism $\alpha^+$ of $\Delta^+$. This gives us a restriction homomorphism $\theta: {\rm Aut}\,\Delta\to{\rm Aut}\,\Delta^+,\;\alpha\mapsto\alpha^+$.

\begin{prop}
The restriction homomorphism $\theta:{\rm Aut}\,\Delta\to{\rm Aut}\,\Delta^+$ is an isomorphism.
\end{prop}

\noindent{\sl Proof.} Suppose that  $\alpha\in\ker\theta$, so that $\alpha$ fixes $\rho_2$ and $\rho_0$. Since $r_1$, acting by conjugation, inverts these two elements, so does its image $r_1\alpha$ under $\alpha$. Thus $r_1(r_1\alpha)^{-1}$ fixes them, so it fixes $\Delta^+$. Since $\Delta^+$ has trivial centraliser in $\Delta$, it follows that  $r_1(r_1\alpha)^{-1}=1$, so $r_1\alpha=r_1$ and hence $\alpha=1$. Thus $\theta$ is a monomorphism.

Each inner automorphism of $\Delta^+$ extends, in the obvious way, to an inner automorphism of $\Delta$, so ${\rm im}\,\theta$ contains ${\rm Inn}\,\Delta^+$. In order to show that $\theta$ is an epimorphism, it is therefore sufficient to show that each of the three outer automorphism classes of $\Delta^+$ corresponding to a generator $X, Y$ or $T$ of $GL_2({\bf Z})$ contains an automorphism of $\Delta^+$ which extends to an automorphism of $\Delta$. In these three cases we can use the automorphisms
\[\alpha_X: \rho_2\mapsto \rho_0^{-1},\,\rho_0\mapsto \rho_2,\quad
\alpha_Y:\rho_2\mapsto\rho_0^{-1},\,\rho_0\mapsto\rho_2\rho_0,\quad
\alpha_T: \rho_2\mapsto\rho_0,\,\rho_0\mapsto\rho_2,\]
which extend to
\[\alpha_x: r_0\mapsto r_2,\, r_1\mapsto r_1,\, r_2\mapsto r_1r_0r_1,\quad
\alpha_y: r_0\mapsto r_2r_1r_2,\, r_1\mapsto r_2, r_2\mapsto r_2r_0r_2,\]
\[\alpha_t: r_0\mapsto r_1r_2r_1,\, r_1\mapsto r_1,\, r_2\mapsto r_1r_0r_1.\]
Thus $\theta$ is an isomorphism. \hfill$\square$

\medskip

It follows from Proposition~2.1 that each automorphism $\beta$ of $\Delta^+$  extends to a unique automorphism $\beta^*$ of $\Delta$, so that $\theta^{-1}$ is given by $\beta\mapsto\beta^*$. Under this isomorphism, ${\rm Inn}\,\Delta^+$ is identified with a characteristic subgroup of ${\rm Inn}\,\Delta$, and thus a normal subgroup of ${\rm Aut}\,\Delta$. By Proposition~2.1, this index $2$ inclusion of ${\rm Inn}\,\Delta^+$ in ${\rm Inn}\,\Delta$ induces an epimorphism $\varepsilon: {\rm Out}\,\Delta^+\to{\rm Out}\,\Delta,\; [\beta]\mapsto[\beta^*]$, corresponding under the isomorphisms ${\rm Out}\,\Delta^+\cong GL_2({\bf Z})$ and ${\rm Out}\,\Delta\cong PGL_2({\bf Z})$ to the natural epimorphism $GL_2({\bf Z})\to PGL_2({\bf Z})$. Its kernel, which has order $2$ and corresponds to the subgroup $\{\pm I\}$ of $GL_2({\bf Z})$, is generated by the outer automorphism class $[\alpha_r]$ of $\Delta^+$; here $\alpha_r$ is the automorphism of $\Delta^+$ defined in Example~2, which extends to an inner automorphism of $\Delta$. To summarise, we have:

\begin{prop}
The extension mapping $\beta\mapsto\beta^*$ induces an epimorphism
\[\varepsilon: {\rm Out}\,\Delta^+\to{\rm Out}\,\Delta,\; [\beta]\mapsto[\beta^*],\]
with $\ker \varepsilon=\langle[\alpha_r]\rangle \cong C_2$. \hfill$\square$
\end{prop}

\begin{cor} If an oriented hypermap is reflexible, then so is its image under any operation on oriented hypermaps.
\end{cor}

\noindent{\sl Proof}. This follows immediately from the fact that $[\alpha_r]$ is in the centre of ${\rm Out}\,\Delta^+$, so that the corresponding reflection operation $\omega_r$ is in the centre of $\Omega^+$. \hfill$\square$

\medskip

In fact, since $[\alpha_r]$ is the only non-identity central element of ${\rm Out}\,\Delta^+$, Corollary~3.2 does not apply to any non-identity operation in $\Omega^+$ other than $\omega_r$. Similarly, since the centre of ${\rm Out}\,\Delta$ is trivial, there is no analogous result for operations on unoriented hypermaps.

\section{Periodic elements of $GL_2({\bf Z})$ and $PGL_2({\bf Z})$.}

The following result concerning periodic elements (those of finite order) is well-known (see
~\cite{Mes}
 for part (a)), but for completeness we will give a proof. 

\begin{prop}
 {\rm(a)} The periodic elements of $GL_2({\bf Z})$ consist of the identity, three conjugacy classes of involutions represented by the matrices
\[-I=\Bigl(\,\begin{matrix}-1&0\\ 0&-1\end{matrix}\,\Bigr),\quad
T=\Bigr(\,\begin{matrix}0&1\\ 1&0\end{matrix}\,\Bigr)\quad {\it and}\quad
XT=\Bigl(\,\begin{matrix}-1&0\\ 0&1\end{matrix}\,\Bigr),\]
and one class each of elements of orders $3, 4$ and $6$ represented by
\[Y^2=\Bigl(\,\begin{matrix}-1&-1\\ 1&0\end{matrix}\,\Bigr),\quad
X=\Bigl(\,\begin{matrix}0&-1\\ 1&0\end{matrix}\,\Bigr)\quad {\it and}\quad
Y=\Bigl(\,\begin{matrix}0&-1\\ 1&1\end{matrix}\,\Bigr).\]

\noindent{\rm(b)}
The periodic elements of $PGL_2({\bf Z})$ consist of the identity, three conjugacy classes of involutions represented by the images $x, t$ and $xt$ of $X, T$ and $XT$, and one conjugacy class of elements of order $3$ represented by the image $y$ of $Y$.
\end{prop}

\noindent{\sl Proof.} (a) The presentation~$(2.1)$ shows that $GL_2({\bf Z})$ is the free product of two dihedral groups $\langle X, T\mid X^4=T^2=(XT)^2=1\rangle$ and $\langle Y, T\mid Y^6=T^2=(YT)^2=1\rangle$ of order $8$ and $12$, amalgamating a Klein four-group generated by $X^2\;(=Y^3=-I)$ and $T$. The torsion theorem for free products with amalgamation
~\cite[Ch.~IV, Theorem~2.7]{LS}
 states that in such a group the periodic elements are the conjugates of those in the factors. In our case, taking account of conjugacy within the two dihedral groups, this implies that the periodic elements of $GL_2({\bf Z})$ are the conjugates of $I,\, -I,\, T,\, XT,\, Y^2,\, X$ and $Y$. (Note that $YT\sim Y^3T=X^2T\sim T$.) These elements (and hence their conjugates) have orders $1,\, 2,\, 2,\, 2,\, 3,\, 4$ and $6$ respectively. The involutions $-I,\, T$ and $XT$ are mapped to distinct elements in the abelianised group $GL_2({\bf Z})^{\rm ab}\cong C_2\times C_2$, so they lie in three different conjugacy classes. The remaining listed elements have different orders, so they also lie in different classes.

(b) Since the natural epimorphism $GL_2({\bf Z})\to PGL_2({\bf Z})$ is finite-to-one, the periodic elements of $PGL_2({\bf Z})$ are the images of those in $GL_2({\bf Z})$. Since conjugate elements of $GL_2({\bf Z})$ have conjugate images in $PGL_2({\bf Z})$, it is sufficient to consider the images of $I,\, -I,\, T,\, XT,\, Y^2,\, X$ and $Y$. The elements $\pm I$ map to the identity, while $T,\, XT$ and $X$ map to the involutions $t,\, xt$ and $x$ in $PGL_2({\bf Z})$. Since these have different images in $PGL_2({\bf Z})^{\rm ab}\cong C_2\times C_2$, they represent three different conjugacy classes. The elements $Y$ and $Y^2$ map to $y$ and $y^2=y^t$, giving a single conjugacy class of elements of order $3$.\hfill$\square$

\medskip

\noindent{\bf Remarks. 1.} One can reverse the process in the proof of (b) to show how the conjugacy classes of periodic elements in $PGL_2({\bf Z})$ lift back to $GL_2({\bf Z})$. The class $\{1\}$ lifts to the union of two classes $\{I\}$ and $\{-I\}$. The classes containing $x$, $t$ and $xt$ each lift to a single class containing $X$, $T$ or $XT$, since $-X=X^{-1}=X^T$, $-T=T^X$ and $-XT=(XT)^X$. The class containing $y$ lifts to the union of the two classes containing $Y$ and $Y^2$, since $-Y=Y^{-2}=(Y^2)^T$.

\smallskip 

\noindent{\bf 2.} One can also prove (b), independently of (a), by using the facts that $PSL_2({\bf Z})$, which has index $2$ in $PGL_2({\bf Z})$, is the free product of $\langle x\rangle\cong C_2$ and $\langle y\rangle\cong C_3$, or that $PGL_2({\bf Z})$ is the free product of $\langle x,t\rangle\cong C_2\times C_2$ and $\langle y,t\rangle\cong S_3$ amalgamating $\langle t\rangle\cong C_2$.

\smallskip 

\noindent{\bf 3.} Similar arguments show that the maximal finite subgroups of $GL_2({\bf Z})$ are the dihedral groups of order $8$ and $12$ conjugate to $\langle X, T\rangle$ and $\langle Y, T\rangle$, while those of $PGL_2({\bf Z})$ are their images, Klein four-groups conjugate to $\langle x, t\rangle$ and dihedral groups of order $6$ conjugate to $\langle y,t\rangle$. 

\medskip

The following result allows easy recognition of periodic elements in $GL_2({\bf Z})$ and $PGL_2({\bf Z})$:

\begin{prop}
{\rm(a)} An element $A\in GL_2({\bf Z})$ is periodic if and only if 
\vskip2pt
\item{\rm(i)} $A= I$ or $-I$ (in which cases $A$ has order $1$ or $2$), or
\vskip2pt
\item{\rm(ii)} $A$ has determinant $1$ and trace $-1, 0$ or $1$ (in which cases $A$ has order $3, 4$ or $6$), or
\vskip2pt
\item{\rm(iii)}  $A$ has determinant $-1$ and trace $0$ (in which case $A$ has order $2$).
\vskip5pt
\item{\rm (b)} An element $a=\pm A\in PGL_2({\bf Z})$ is periodic if and only if 
\vskip2pt
\item{\rm(i)} $A=\pm I$ (in which case $a$ has order $1$), or
\vskip2pt
\item{\rm(ii)} $A$ has determinant $1$ and trace $0$ or $\pm 1$ (in which cases $a$ has order $2$ or $3$), or
\vskip2pt
\item{\rm(iii)}  $A$ has determinant $-1$ and trace $0$ (in which case $a$ has order $2$).
\end{prop}

\noindent{\sl Proof.} (a) By Proposition~3.1(a), any periodic element $A$ is conjugate to $\pm I$ (so that $A=\pm I$) or to one of $T, XT, Y^2, X$ or $Y$. Since determinant, trace and order are invariant under conjugacy, it follows that in the last five cases $A$ satisfies (ii) or (iii). Conversely, if $A$ satisfies (ii) or (iii) then by solving the characteristic equation one sees that the eigenvalues of $A$ are a pair of distinct roots of unity of the stated order, so diagonalising $A$ (over $\bf C$) shows that $A$ also has this order.

\smallskip

(b) The proof for $PGL_2({\bf Z})$ is similar, using Proposition~3.1(b). \hfill$\square$

\medskip

This proof shows that in cases (ii) and (iii) of Proposition~3.2(a), the determinant and trace of $A$ determine its conjugacy class uniquely, except when $\det A=-1$ and ${\rm tr}\, A=0$, since in this case $A$ could be conjugate to $T$ or $XT$. A similar phenomenon occurs in $PGL_2({\bf Z})$, where the conjugacy classes containing $t$ and $xt$ cannot be distinguished by determinant and trace. In either case, the simplest way to determine the conjugacy class of an element is to reduce its entries mod~$(2)$. The kernel of the natural epimorphism $GL_2({\bf Z})\to GL_2({\bf Z}/2{\bf Z})$, the principal congruence subgroup of level $2$, contains $XT$ but not $T$, so $A={a\; b\choose c\;d}$ is conjugate to $XT$ if $b$ and $c$ are even (so $a$ and $d$ are odd), and otherwise $A$ is conjugate to $T$. The same criterion determines whether $\pm A\in PGL_2({\bf Z})$ is conjugate to $xt$ or $t$.

\section{Periodic operations on oriented hypermaps}

Using the isomorphism $GL_2({\bf Z})\cong {\rm Out}\,\Delta^+=\Omega^+$ one can interpret the elements of $GL_2({\bf Z})$ as operations on oriented hypermaps. We will consider the operation corresponding to a representative of each of the conjugacy classes listed in Proposition~3.1(a).  Unfortunately, the extension ${\rm Aut}\,\Delta^+$ of ${\rm Inn}\,\Delta^+$ by ${\rm Out}\,\Delta^+$ is not split, so there is no single subgroup of ${\rm Aut}\,\Delta^+$  which provides representatives of the required outer automorphism classes: instead, one has to make arbitrary choices of representatives.

The automorphisms
\[\alpha_X: \rho_2\mapsto \rho_0^{-1},\,\rho_0\mapsto \rho_2\quad\hbox{\rm and}\quad 
\alpha_T: \rho_2\mapsto\rho_0,\,\rho_0\mapsto\rho_2\]
used in the proof of Proposition~2.1 satisfy the relations
\[\alpha_X^4=\alpha_T^2=(\alpha_X\alpha_T)^2=1,\]
and generate a dihedral group of order $8$ in ${\rm Aut}\,\Delta^+$ which is mapped isomorphically into ${\rm Out}\,\Delta^+$. We can therefore use the elements of this group as representatives of the corresponding outer automorphism classes, and in particular those classes corresponding to the periodic elements $X, -I, T$ and $XT$ in Proposition~3.1(a).

As in the proof of Proposition~2.1, the outer automorphism class corresponding to the matrix $X$ is represented by the automorphism $\alpha_X: \rho_2\mapsto\rho_0^{-1},\, \rho_0\mapsto \rho_2$ of order $4$; the corresponding operation $\omega_X\in\Omega^+$ transposes hypervertices and hyperfaces, reversing the orientation around the latter. This operation is conjugate (under $\omega_T$) to its inverse, which also transposes hypervertices and hyperfaces but reverses the orientation around the former.

The class corresponding to $-I$ is represented by the automorphism $\alpha_{-I}=\alpha_X^2=\alpha_r$ of $\Delta^+$ inverting $\rho_2$ and $\rho_0$ (see Example~2). This is the restriction to $\Delta^+$ of the inner automorphism $\iota_{r_1}$ of $\Delta$ induced by $r_1$, and the corresponding operation $\omega_{-I}=\omega_X^2=\omega_r$ on oriented hypermaps simply reverses their orientation.

As in the proof of Proposition~2.1, the class corresponding to $T$ is represented by the automorphism $\alpha_T$ transposing $\alpha_2$ and $\alpha_0$. The corresponding operation $\omega_T\in\Omega^+$ is duality, transposing hypervertices and hyperfaces while preserving hyperedges and orientation; among its conjugates are the other two dualities, transposing hyperedges with hypervertices or hyperfaces, and the operation $\omega_{-T}$ which transposes hypervertices and hyperfaces while reversing orientation (see \S 6 for a further discussion of duality).

The class corresponding to $XT$ is represented by the automorphism $\alpha_{XT}=\alpha_X\alpha_T$ inverting $\rho_2$ and fixing $\rho_0$; this corresponds to a `twisting' operation $\omega_{XT}=\omega_X\omega_T$ which reverses the orientation around the hyperfaces but preserves it around the hypervertices. This operation is one of six conjugate twisting operations, each of which reverses the orientation around hypervertices, hyperedges or hyperfaces, while preserving it around one other of these three constituents of a hypermap. These twisting operations can change the valencies of the third constituent, and therefore (like $\omega_X$) they can change the type and the genus of a hypermap, whereas these are all preserved by $\omega_{-I}$ and $\omega_T$.

As in the preceding examples, the matrix $Y^2$ of order $3$ can also be represented by an automorphism of the same order, namely
\[\alpha_{Y^2}: \rho_2\mapsto\rho_0^{-1}\rho_2^{-1}=\rho_1,\; \rho_0\mapsto\rho_2.\]
This permutes $\rho_2,\, \rho_1$ and $\rho_0$ in a $3$-cycle, thus inducing the triality operation $\omega_{Y^2}$ which permutes hyperfaces, hyperedges and hypervertices in a $3$-cycle; this operation is conjugate to its inverse.

In the case of the matrix $Y$ of order $6$, however, we cannot find a corresponding automorphism of the same order since, as shown by Meskin~\cite{Mes}, $F_2$ has no automorphisms of order~$6$. Instead, we must use an automorphism $\alpha_Y\in{\rm Aut}\,\Delta^+$ of infinite order such that $[\alpha_Y]$ has order $6$ in ${\rm Out}\,\Delta^+$: for instance, we can choose the automorphism
\[\alpha_Y:\rho_2\mapsto \rho_0^{-1},\, \rho_0\mapsto \rho_2\rho_0=\rho_1^{-1}\]
appearing in the proof of Proposition~2.1, so that $\alpha_Y^2=\alpha_{Y^2}\circ\iota_{\rho_1}^{-1}$, $\alpha_Y^3=\alpha_{-I}\circ\iota_{\rho_1}^{-1}$ and $\alpha_Y^6=\iota_{\rho_1}^{-2}$ where $\iota_{\rho_1}$ is the inner automorphism induced by $\rho_1$. The corresponding operation $\omega_Y$ sends hyperfaces and hypervertices to hypervertices and hyperedges, reversing the orientation around each. Its powers $\omega_Y^2, \omega_Y^3$ and $\omega_Y^6$ induce the operations $\omega_{Y^2}, \omega_{-I}$ and the identity operation $\omega_I$.

\medskip

\noindent{\bf Example 4.} The icosahedron, an orientably regular hypermap of type $(5,2,3)$, can be regarded as a permutation representation $\Delta^+\to G=A_5$ such as
\[\rho_0\mapsto (1,2,3,4,5), \; \rho_1\mapsto (1,2)(3,4),\; \rho_2\mapsto(2,5,4).\]
Applying $\omega_{XT}$, we have $\rho_0\mapsto (1,2,3,4,5)$ again but now $\rho_2\mapsto(2,5,4)^{-1}=(2,4,5)$, so $\rho_1=\rho_0^{-1}\rho_2^{-1}\mapsto (1,4,3,5,2)$, giving an orientably regular hypermap of type $(5,5,3)$. Repeated use of this type of calculation determines the orbits ${\cal O}_1$ and ${\cal O}_2$ in Example~3.

\medskip

This section has shown how the periodic elements of $GL_2({\bf Z})$, such as $X, Y$ and $T$, induce operations on oriented hypermaps. Since these three elements generate $GL_2({\bf Z})$, this allows one to find the operation $\omega_A$ corresponding to {\sl any\/} $A\in GL_2({\bf Z})$: first one can use elementary row operations to express $A$ as a word in the generators $X, Y$ and $T$ (see \S 8), and then by replacing these with $\alpha_X, \alpha_Y$ and $\alpha_T$, or with $\omega_X, \omega_Y$ and $\omega_T$, one can find an automorphism $\alpha_A$ of $F_2$ acting as $A$ on $F_2^{\rm ab}$, together with the corresponding operation $\omega_A$ on oriented hypermaps. Of course, $\alpha_A$ is not uniquely determined by $A$, but any two such automorphisms differ by an inner automorphism of $F_2$, so $\omega_A$ is unique.

\section{Periodic operations on all hypermaps}

In a similar way, the isomorphism $PGL_2({\bf Z})\cong {\rm Out}\,\Delta=\Omega$ and Proposition~2.1(b) can be used to interpret the periodic elements of $PGL_2({\bf Z})$ as operations on all hypermaps, ignoring orientation if it exists. These operations are induced by automorphisms of $\Delta$ which are the extensions of the automorphisms of $\Delta^+$ used in the preceding section.

The involution $x$ corresponds to the outer automorphism class of $\Delta$ containing the automorphism
\[\alpha_x: r_0\mapsto r_2,\; r_1\mapsto r_1,\; r_2\mapsto r_1r_0r_1\]
used in the proof of Proposition~2.1. This restricts to the automorphism $\alpha_X$ of $\Delta^+$, and the corresponding operation $\omega_x$ acts on all hypermaps as the operation $\omega_X$ described in \S 4, transposing hypervertices and hyperfaces, and reversing orientation around the latter. Then $\alpha_x^2$ is the inner automorphism $\iota_{r_1}$, which induces the trivial operation on unoriented hypermaps though it acts on oriented hypermaps as orientation-reversal $\omega_{-I}$.

The involution $t$ is represented by the automorphism
\[\alpha_t: r_0\mapsto r_1r_2r_1,\; r_1\mapsto r_1,\; r_2\mapsto r_1r_0r_1\]
of $\Delta$, which restricts to $\alpha_T$ on $\Delta^+$; this is the composition of transposing $r_0$ and $r_2$ with conjugation by $r_1$, so the corresponding operation $\omega_t$ is the duality operation $\omega_{02}$ transposing hyperfaces and hypervertices (see Example~1).

The involution $xt$ is represented by the automorphism
\[\alpha_{xt}: r_0\mapsto r_1r_0r_1,\; r_1\mapsto r_1,\; r_2\mapsto r_2,\]
which restricts to $\alpha_{XT}$ on $\Delta^+$; as before, this corresponds to a twisting operation $\omega_{xt}$ which preserves the hypervertices and reverses the orientation around each hyperface.

The element $y$ is represented by the automorphism
\[\alpha_y: r_0\mapsto r_2r_1r_2,\; r_1\mapsto r_2,\; r_2\mapsto r_2r_0r_2.\]
which restricts to $\alpha_Y$ on $\Delta^+$. This is the composition of the $3$-cycle $(r_0, r_1, r_2)$ with $\iota_{r_2}$, so the corresponding operation $\omega_y$ acts on isomorphism classes of hypermaps by permuting hypervertices, hyperedges and hyperfaces in a $3$-cycle (see Example~1).

The remarks at the end of \S 4 also apply here, with $x, y$ and $t$ replacing $X, Y$ and $T$, allowing one to determine an automorphism $\alpha_a$ and the operation $\omega_a$ corresponding to any element $a=\pm A\in PGL_2({\bf Z})$.

\section{Duality operations}

The classical duality for maps transposes vertices and faces, and its most obvious extension to oriented hypermaps is the operation $\omega_T$, induced by the automorphism
\[\alpha_T: \rho_2\mapsto\rho_0,\; \rho_0\mapsto\rho_2\]
of $\Delta^+$ which sends $\rho_1=\rho_0^{-1}\rho_2^{-1}$ to $\rho_2^{-1}\rho_0^{-1}=\rho_2^{-1}\rho_1\rho_2$. This operation is conjugate in $\Omega^+$, under $\omega_Y$ and $\omega_Y^2$, to the other two dualities of oriented hypermaps which transpose hyperedges with hypervertices or hyperfaces by transposing $\rho_1$ with $\rho_0$ or $\rho_2$.

Similarly, in the case of unoriented hypermaps, the natural definition of duality is the operation $\omega_{02}$ induced by the automorphism
\[\alpha_{02}:r_0\mapsto r_2,\; r_1\mapsto r_1,\; r_2\mapsto r_0\]
of $\Delta$ (see Example~1). This is conjugate in $\Omega$ to two other duality operations $\omega_{01}$ and $\omega_{12}$ which transpose hyperedges with hypervertices or hyperfaces. However the restriction $\alpha^+_{02}$ of $\alpha_{02}$ to $\Delta^+$ is not $\alpha_T$ but rather the automorphism
\[\alpha_{-T}=\alpha_T\circ\alpha_{-I}=\alpha_{-I}\circ\alpha_T: \rho_2\to \rho_0^{-1},\; \rho_0\to\rho_2^{-1},\]
corresponding to the matrix $-T$, which inverts $\rho_1$. This induces on oriented hypermaps the operation $\omega_{-T}=\omega_T\circ\omega_{-I}=\omega_{-I}\circ\omega_T$ which combines the classical duality $\omega_T$ with a reversal of orientation. Similarly the duality operations $\omega_{01}$ and $\omega_{02}$, when restricted to oriented hypermaps, reverse the orientation. (This problem does not arise in connection with the two triality operations corresponding to the $3$-cycles in $S_3$, since their restrictions preserve orientation.)

The two duality operations $\omega_T$ and $\omega_{-T}$ on oriented hypermaps are distinct, but since they are conjugate in $\Omega^+$ (under $\omega_X$, since $T^X=-T$) many of their properties are similar: for instance, they give rise to the same generalised chirality groups, as defined in \S 7. The same applies to the two other pairs of duality operations on oriented hypermaps, transposing hyperedges with hypervertices or hyperfaces: in each case the two operations differ by an orientation-reversal, and are conjugate to each other in $\Omega^+$.

\section{Generalised chirality groups}

Let $\omega$ be an operation of order $2$ on hypermaps or on oriented hypermaps, induced by an automorphism $\alpha$ of a group $\Gamma=\Delta$ or $\Delta^+$ respectively, so that $\alpha$ is not inner but $\alpha^2$ is, and let $\cal H$ be a regular or orientably regular hypermap corresponding to a normal subgroup $H$ of $\Gamma$. Then ${\cal H}^{\omega}$ is the hypermap corresponding to the normal subgroup $H^{\alpha}$, and since $\alpha$ transposes $H$ and $H^{\alpha}$ and preserves the normal subgroups $HH^{\alpha}$ and $H\cap H^{\alpha}$, the quotient groups $HH^{\alpha}/H$, $HH^{\alpha}/H^{\alpha}$, $H/H\cap H^{\alpha}$ and $H^{\alpha}/H\cap H^{\alpha}$ are all isomorphic to each other. We call this common quotient group the {\sl generalised chirality group} $\chi_{\omega}({\cal H})$ of $\cal H$ corresponding to $\omega$; it is independent of the choice of an automorphism $\alpha$ representing $\omega$ since any two of them differ by an inner automorphism of $\Gamma$. The order $|\chi_{\omega}({\cal H})|$ of this group is the {\sl generalised chirality index} of $\cal H$. Since $HH^{\alpha}$ and $H\cap H^{\alpha}$ correspond to the largest $\omega$-invariant regular quotient of $\cal H$ and the smallest $\omega$-invariant regular covering of $\cal H$, both the generalised chirality group and the index are measures of the extent to which $\cal H$ fails to be $\omega$-invariant. In the cases where $\omega$ is the chirality or duality operation $\omega_{-I}$ or $\omega_T$ on oriented hypermaps, $\chi_{\omega}({\cal H})$ is respectively the chirality group studied in~\cite{BJNS} or the duality group studied in~\cite{Pin}.

If $\omega$ and $\omega'$ are operations induced by automorphisms $\alpha$ and $\alpha'$ of $\Gamma$, and are conjugate, say $\omega'=\tau\omega\tau^{-1}$ for some operation $\tau$ induced by an automorphism $\beta\in{\rm Aut}\,\Gamma$, then 
\[\chi_{\omega'}({\cal H})\cong HH^{\alpha'}/H=HH^{\beta\alpha\beta^{-1}}/H\cong H^{\beta}H^{\beta\alpha}/H^{\beta}\cong\chi_{\omega}({\cal H}^{\tau}).\]
Thus the groups which can arise as generalised chirality groups for two conjugate operations are the same, so this set of groups depends only on the conjugacy class of an operation. In particular, the three duality operations $\omega_{ij}$ on hypermaps are conjugate to each other in $\Omega$, so they give rise to the same generalised chirality groups. The same applies to the six duality operations on oriented hypermaps (two for each pair $i, j$, preserving or reversing orientation), since they are all conjugate to each other in $\Omega^+$, as shown in \S 6. For instance, it is shown in~\cite{BJNS} that every finite abelian group can appear as the chirality group $\chi_{-I}({\cal H})$ of a finite hypermap $\cal H$; a similar argument in~\cite{Pin} gives the corresponding result for the duality operation $\omega_T$, and hence (by conjugacy) for all six duality operations on oriented hypermaps.
 
The notion of a generalised chirality group or index does not extend so easily to the case where $\omega$ has order greater than $2$. The natural analogues of $HH^{\alpha}$ and $H\cap H^{\alpha}$ would then be the group $H^{\langle\alpha\rangle}$ generated by all the images $H^{\alpha^i}$ of $H$ under powers of $\alpha$, and the group $H_{\langle\alpha\rangle}$ which is the intersection of those images. These correspond to the largest $\omega$-invariant quotient of $\cal H$ and the smallest $\omega$-invariant covering of $\cal H$. In general, the {\sl upper\/} and {\sl lower generalised chirality groups\/} $\chi^{\omega}({\cal H})=H^{\langle\alpha\rangle}/H$ and $\chi_{\omega}({\cal H})=H/H_{\langle\alpha\rangle}$ need not be isomorphic. For instance, for each prime $p\geq 5$ the group $L:=L_2(p)=PSL_2({\bf Z}/p{\bf Z})$ is a simple quotient of the modular group $PSL_2({\bf Z})$ (by its principal congruence subgroup of level $p$), and hence of $\Delta^+$. The corresponding kernel $H$ is normal not only in $\Delta^+$ but also in $\Delta$ since $L$ has an automorphism inverting both of its generators, so the corresponding hypermap $\cal H$ is reflexible. (In fact, ${\rm Aut}\,{\cal H}\cong L\times C_2$ or $PGL_2(p)$ as $p\equiv 1$ or $p\equiv -1$ mod~$(4)$.) If we take $\omega$ to be the triality operation $\omega_y$ of order $3$ then since $\cal H$ has type $(2,3,p)$ the three hypermaps $\cal H$, ${\cal H}^{\omega}$ and ${\cal H}^{\omega^2}$ are mutually non-isomorphic. It follows from the simplicity of $L$ that $HH^{\alpha}=\Delta^+$ and $\Delta^+/H\cap H^{\alpha}\cong L\times L$. Since $L$ is non-abelian and simple the two direct factors of $L\times L$ are its only normal subgroups with quotients isomorphic to $L$, so $H^{\alpha^2}\not\geq H\cap H^{\alpha}$. Thus $\Delta^+/H_{\langle\omega\rangle} \cong L\times L\times L$ and so $\chi_{\omega}({\cal H})\cong L\times L$ whereas $\chi^{\omega}({\cal H})\cong \Delta^+/H \cong L$. In this case, $|\chi_{\omega}({\cal H})|>|\chi^{\omega}({\cal H})|$, but on the other hand, if we take $\cal K$ to be the hypermap corresponding to $H\cap H^{\alpha}$ then $\chi_{\omega}({\cal K})\cong L$ and $\chi^{\omega}({\cal K})\cong L\times L$, so the inequality is reversed.

\section{Canonical forms}

Although we are mainly concerned here with operations of finite order, the majority of operations in  $\Omega$ or $\Omega^+$ have infinite order: in each case, for instance, they form infinitely many conjugacy classes, whereas Proposition~3.1 shows that there are only finitely many conjugacy classes of operations of finite order. However, each of these groups is generated by operations of finite order, for instance $\omega_x, \omega_y, \omega_t$ and $\omega_X, \omega_Y, \omega_T$ respectively. We can use the free product structure of $PSL_2({\bf Z})$ to give a unique expression, or canonical form, for each element of $\Omega$ and $\Omega^+$ in terms of these generators.

The group $PSL_2({\bf Z})$ is the free product of $\langle x\rangle\cong C_2$ and $\langle y\rangle\cong C_3$, so each element $a\in PSL_2({\bf Z})$ can be expressed as a unique reduced word
\[a=w(x,y)=x^{\delta_1}y^{\varepsilon_1}\ldots\, x^{\delta_k}y^{\varepsilon_k}\eqno(8.1)\]
in $x$ and $y$, where `reduced' means that there is no internal cancellation, so in this case $\delta_1=0$ or $1$, $\delta_2=\ldots=\delta_k=1$, each $\varepsilon_1,\ldots, \varepsilon_{k-1}=1$ or $-1$, and $\varepsilon_k=0, 1$ or $-1$. Since $PGL_2({\bf Z})$ is a semidirect product of $PSL_2({\bf Z})$ by $\langle t\rangle\cong C_2$, each element $a\in PGL_2({\bf Z)}\setminus PSL_2({\bf Z})$ has the form $a=w(x,y)t$ for a unique reduced word $w(x,y)$. Thus each element $a\in PGL_2({\bf Z})$ has the unique form
\[a=w(x,y)t^{\eta}\eqno(8.2)\]
where $w(x,y)$ is a reduced word and $\eta=0$ or $1$; by the isomorphism $a\mapsto\omega_a$ between $PGL_2({\bf Z})$ and $\Omega$, the same is true for each operation $\omega_a\in\Omega$, with $\omega_x, \omega_y$ and $\omega_t$ replacing $x, y$ and $t$. 

Similarly, if $A\in GL_2({\bf Z})$ then its image $a\in PGL_2({\bf Z})$ under the natural epimorphism $GL_2({\bf Z})\to PGL_2({\bf Z})$ has a canonical form $(8.2)$, so $A$ has a canonical form
\[A=\pm w(X,Y)T^{\eta}\eqno(8.3)\]
where again $\eta=0$ or $1$ and $w$ is a reduced word in the sense defined above, with the same restrictions on the exponents $\delta_i$ and $\varepsilon_i$. This gives a corresponding canonical form for each operation $\omega_A\in\Omega^+$.

There is a simple algorithm for putting any element of $GL_2({\bf Z})$ or $PGL_2({\bf Z})$, and hence any operation in $\Omega^+$ or $\Omega$, into canonical form. If $A\in GL_2({\bf Z})$ then multiplying $A$ on the left by suitable powers of the matrices
\[\Bigl(\,\begin{matrix}1&1\\ 0&1\end{matrix}\,\Bigr)=-XY
\quad{\rm and}\quad
\Bigl(\,\begin{matrix}1&0\\ 1&1\end{matrix}\,\Bigr)=XY^{-1}\]
corresponds to using row operations to reduce $A$ to $\pm X^{\delta}T^{\eta}$ for some $\delta, \eta \in\{0, 1\}$. This yields an equation of the form $W(X,Y)A=\pm T^{\eta}$ for some word $W$ in $X$ and $Y$, with $\eta=0$ or $1$, so $A=\pm W(X,Y)^{-1}T^{\eta}$, and a finite process converts this expression to the form $(8.3)$ where $w$ is a reduced word. Essentially the same algorithm also puts any element $a=\pm A\in PGL_2({\bf Z})$ into its canonical form $(8.2)$.

\medskip

\noindent{\bf Example 5.} Let
\[A=\Bigl(\,\begin{matrix}-2&-3\\ 1&2\end{matrix}\,\Bigr)\in GL_2({\bf Z}).\]
The reduction process
\[\Bigl(\,\begin{matrix}1&1\\ 0&1\end{matrix}\,\Bigr)
\Bigl(\,\begin{matrix}-2&-3\\1&2\end{matrix}\,\Bigr)
=\Bigl(\,\begin{matrix}-1&-1\\1&2\end{matrix}\,\Bigr),
\quad
\Bigl(\,\begin{matrix}1&0\\ 1&1\end{matrix}\,\Bigr)
\Bigl(\,\begin{matrix}-1&-1\\ 1&2\end{matrix}\,\Bigr)
=\Bigl(\,\begin{matrix}-1&-1\\ 0&1\end{matrix}\,\Bigr),\]
\[\Bigl(\,\begin{matrix}1&1\\ 0&1\end{matrix}\,\Bigr)
\Bigl(\,\begin{matrix}-1&-1\\ 0&1\end{matrix}\,\Bigr)
=\Bigl(\,\begin{matrix}-1&0\\ 0&1\end{matrix}\,\Bigr)\]
gives
\[(-XY)(XY^{-1})(-XY)A=XT.\]
Hence
\begin{eqnarray*}
A & = & (-XY)^{-1}(XY^{-1})^{-1}(-XY)^{-1}XT \\
& = & Y^{-1}X^{-1}YX^{-1}Y^{-1}X^{-1}XT \\
& = & Y^{-1}XYXY^{-1}T
\end{eqnarray*}
in canonical form $(8.3)$, where we have twice replaced $X^{-1}$ with $-X$. Replacing the generators $X, Y$ and $T$ with their corresponding operations gives the canonical form for the operation $\omega_A$.

\medskip

One can characterise the operations of finite order in $\Omega$ and in $\Omega^+$, and determine their conjugacy classes, in terms of these canonical forms, though it is notationally simpler to work with $PGL_2({\bf Z})$ and $GL_2({\bf Z})$. If $a\in PSL_2({\bf Z})$, so that $\eta=0$ and $a=w(x,y)$ in $(8.2)$, then the torsion and conjugacy theorems for free products~\cite[\S IV.1]{LS}, imply that $a$ has finite order if and only if $w$ can be cyclically reduced (by successively cancelling mutually inverse first and last terms) to a power of $x$ or $y$, in which case $a$ is conjugate in $PSL_2({\bf Z})$ to that power. This is equivalent to $w$ being the empty word (so that $a=1$), or to $w$ having the form $w_1(x,y)z\overline{w}_1(x^{-1},y^{-1})$ where $z=x$, $y$ or $y^{-1}$, with $w_1$ a reduced word not ending in a power of $z$, and $\overline{w}_1$ the reversed word of $w_1$. Of course, the conjugacy classes of $y$ and $y^{-1}$ in $PSL_2({\bf Z})$ form a single conjugacy class in $PGL_2({\bf Z})$ since $y^t=y^{-1}$.

An element $a=w(x,y)t\in PGL_2({\bf Z})\setminus PSL_2({\bf Z})$ has finite order if and only if the element $a^2\in PSL_2({\bf Z})$ has finite order. Now $a^2=(wt)^2=w^tw=w(x,y^{-1})w(x,y)$, with even and zero exponent sums in $x$ and $y$, so this cannot be cyclically reduced to $x, y$ or $y^{-1}$. It follows that $a$ has finite order if and only if $a^2$ can be cyclically reduced to $1$, that is, $a^2=1$.
This equivalent to $w(x,y^{-1})=w(x,y)^{-1}=\overline w(x^{-1},y^{-1})=\overline w(x,y^{-1})$ where $\overline w$ is the reverse word of $w$, and by the normal form theorem for free products~\cite[\S IV.1]{LS},  this happens if and only if $w=\overline w$, that is, $w$ is a palindrome. By conjugating $a$ with suitable powers of $x$ and $y$, and then putting the resulting elements in the canonical form $w't$ with $w'\in PSL_2({\bf Z})$, one can reduce the wordlength to the unique minimal cases where $w'=1$ or $x$: for instance, conjugating $a = yt$ with $y^{-1}$ gives $w't=y.yt.y^{-1} = y^2.yt = t$. It follows that $a$ is conjugate to $xt$ if the palindrome $w$ has odd length with middle letter $x$, and otherwise $a$ is conjugate to $t$. 

In the case of $GL_2({\bf Z})$, an element $A$ with canonical form $(8.3)$ has finite order if and only if its image $a=\pm A$ in $PGL_2({\bf Z})$, with canonical form $(8.2)$, has finite order. The above criteria can be therefore applied to $A$ by first applying them to $a$ and then using Remark~1 of \S 3 to deduce the conjugacy class of $A$ from that of $a$. The only possible ambiguity is when $a$ is conjugate to $y$; in this case $A$ is conjugate to $Y$ or $Y^2$ as $A$ has canonical form $w(X,Y)$ or $-w(X,Y)$ in $(8.3)$, since these elements have order $6$ and $3$ respectively.

\medskip

\noindent{\bf Example 6.} If $A$ is as in Example~5 then conjugating by $y^{-1}$, then $x$, and then $y^{-1}$ gives
\[a=y^{-1}xyxy^{-1}t\sim xyxy^{-1}ty^{-1}=xyxy^{-1}yt=xyxt\sim yxtx=yx^2t=yt\sim y^3t=t,\]
so $A$ is conjugate to $T$ (see also Proposition~3.2 and the remarks following it) and hence $\omega_A$ is conjugate to the duality operation $\omega_T$.

\end{document}